\documentclass{article}

\usepackage[latin1]{inputenc}
\usepackage{graphicx}
\usepackage{citesort}
\usepackage{amssymb,amsmath}

\usepackage{epsfig}
\usepackage{algorithm}
\usepackage{algorithmicx}
\usepackage{algpseudocode}
\usepackage{authblk}
\usepackage{fullpage}

\newcommand{\Intel}{Intel\textsuperscript{\textregistered}\ }
\newcommand{\Xeon}{Xeon\textsuperscript{\textregistered}\ }

\newcommand{\NI}{\noindent}
\newcommand{\VS}{{\vspace*{0.05in}\NI}}

\newcommand{\RFEAST}{{\scriptsize R-FEAST}}
\newcommand{\BFEAST}{{\scriptsize Bi-FEAST}}
\newcommand{\pk}{\phantom{kk}}

\newcommand{\complex}{\mathbb{C}}

\newcommand{\domain}{{\cal C}}
\newcommand{\boundary}{\partial {\cal C}}
\newcommand{\ConInt}[1]{\frac{1}{2\pi\eye}\oint_{\boundary}\,{#1}\,dz}

\newcommand{\vv}[1]{\mathbf{#1}}
\newcommand{\vvu}{\vv{u}}
\newcommand{\vvv}{\vv{v}}

\newcommand{\vvx}{\vv{x}}
\newcommand{\vvy}{\vv{y}}

\newcommand{\vve}{\vv{e}}
\newcommand{\vvf}{\vv{f}}

\newcommand{\vvk}{\vv{k}}

\newcommand{\Enorm}[1]{\| {#1} \|_2}
\newcommand{\inverse}[1]{{#1}^{-1}}

\newcommand{\iterate}[2]{#1_{({#2})}}
\newcommand{\eig}{{\rm eig}}
\newcommand{\eps}{\varepsilon}

\newcommand{\TSec}[2]{#1_{#2}}
\newcommand{\BSec}[2]{#1_{{#2}'}}
\newcommand{\XTop}[1]{\TSec{X}{#1}}

\newcommand{\XBot}[1]{\BSec{X}{#1}}

\newcommand{\YTop}[1]{\TSec{Y}{#1}}

\newcommand{\YBot}[1]{\BSec{Y}{#1}}

\newcommand{\GammaTop}[1]{\TSec{\Gamma}{#1}}
\newcommand{\GammaBot}[1]{\BSec{\Gamma}{#1}}
\newcommand{\LambdaTop}[1]{\TSec{\Lambda}{#1}}
\newcommand{\LambdaBot}[1]{\BSec{\Lambda}{#1}}

\newcommand{\ctrans}[1]{{#1}^{H}}
\newcommand{\invctrans}[1]{{#1}^{-H}}

\newcommand{\diag}{{\rm diag}}
\newcommand{\rank}{{\rm rank}}
\newcommand{\spann}{{\rm span}}

\newcommand{\quadf}{\rho}
\newcommand{\reference}{{\rm ref}}

\newcommand{\surM}{\inverse{B}A}

\newcommand{\Rprojector}{{X_\domain}\ctrans{(Y_\domain)} B}
\newcommand{\Lprojector}{{Y_\domain}\ctrans{(X_\domain)} \ctrans{B}}
\newcommand{\AI}{\eta}

\newcommand{\Span}{{\rm span}}

\newcommand{\eye}{{\iota}}

\newcommand{\bydef}{\stackrel{\mathrm{def}}{=}}

\begin{document}


\title{A New Highly Parallel Non-Hermitian Eigensolver}

\author[1]{Ping Tak Peter Tang}
\author[2]{James Kestyn}
\author[2]{Eric Polizzi}
\affil[1]{Intel Corporation, 2200 Mission College Blvd, Santa Clara, CA 95054}
\affil[2]{Department of Electrical and Computer Engineering, University of Massachusetts, Amherst, 
          MA 01003}

\maketitle


\newcommand{\AMSsub}[1]{\noindent\textbf{AMS Classifications: }#1}

\begin{abstract}

Calculating portions of eigenvalues and eigenvectors of matrices or matrix pencils has many applications. An
approach to this calculation for Hermitian problems based on a density matrix has been proposed in 2009
and a software package called FEAST has been developed. The density-matrix approach allows FEAST's implementation
to exploit a key strength of modern computer architectures, namely, multiple levels of parallelism.
Consequently, the software package has been well received and subsequently commercialized. 
A detailed theoretical analysis of Hermitian FEAST has also been established very recently. This paper
generalizes the FEAST algorithm and theory, for the first time, to tackle non-Hermitian problems. Fundamentally,
the new algorithm is basic subspace iteration or Bauer bi-iteration, except applied with a novel accelerator
based on Cauchy integrals. The resulting algorithm retains the multi-level parallelism of Hermitian
FEAST, making it a valuable new tool for large-scale computational science and engineering problems on
leading-edge computing platforms.


\end{abstract}


\section{Introduction}
\label{sec:intro}

Generalized non-Hermitian eigenvalue problems of the form
$A x = \lambda B x$ arise in many important applications of applied
sciences and engineering that include economic modeling, Markov chain modeling, structural
engineering, fluid mechanics, material science, and more 
(see~\cite{saad-chebyshev-1984,saad-eigenvalue-problems-2011} for example).
Solving complex-symmetric (still non-Hermitian) eigenvalue problems are crucial in 
modeling open systems based on the perfectly matched layer (PML) technique that is
staple tool in electromagnetics
nanoelectronics~\cite{odermatt-luisier-witzigmann-2005}, and micro electromechanical
systems MEMS~\cite{bindel-govindjee-2005}.
As a tool in numerical linear algebra, non-Hermitian eigensolvers
are kernels to non-linear eigenvalue problems
such as quadratic or polynomial eigenvalue 
problems~\cite{saad-eigenvalue-problems-2011,tisseur-meerbergen-2001}.
More generally, advances in high-performance and big-data computing 
will only increase
the use for general eigenvalue solvers in areas such as bioinformatics, social network,
data mining, just to name a few.
Compared to the Hermitian case, the arsenal of solvers available for non-Hermitian
eigenproblems are much more meager.\footnote{
See {\scriptsize\tt www.netlib.org/utk/people/JackDongarra/la-sw.htm}}
Any addition to the software toolbox for the general scientific computing is 
therefore always timely and welcome.

\VS

For eigenproblems of moderate size, robust solvers are well developed and widely available~\cite{LAPACK-1999} and
are sometimes referred to as direct solvers~\cite{demmel-numerical-linear-algebra}. These solvers
typically calculate the entire spectrum of the matrix or matrix pencil in question. In
many applications, especially for those where the underlying linear systems are large and sparse,
often only selected regions of the spectrum are of interest.
A new approach for
these calculations for Hermitian matrices and matrix pencils based on density
matrices has been proposed recently~\cite{polizzi-2009}. 
Unlike well-known Krylov subspace methods
(see for example~\cite{bai-etal-template-2000,cullum-willoughby-1985,lehoucq-sorensen-1996,parlett-1998})
which maintain subspaces of increasing dimensions,
the FEAST algorithm maintains a basis for a fixed-dimension subspace but updates it
per iteration. In this view, it is similar to the non-expanding subspace version of an eigensolver
based on trace minimization~\cite{sameh-tong-2000,sameh-wisniewski-1982} but with a different subspace update
strategy.  From an implementation point of view, this new approach is similar to spectral
divide-and-conquer~\cite{bai-etal-1997,bai-demmel-1998} in that the calculation is expressed in
terms of high-level building blocks that can much better exploit the advantages of modern computing
architectures. In this case, the high-level building block is a numerical-quadrature based technique
to approximate an exact spectral projector. This building block consists of solving independent
linear systems, each for multiple right hand sides. A software package 
FEAST\footnote{Available at {\scriptsize\tt www.ecs.umass.edu/\~{}polizzi/feast}.}
based on this approach has been made available since 2009.
A comprehensive theoretical analysis of Hermitian FEAST has been completed very 
recently~\cite{tang-polizzi-2013-FEAST} by two of the authors of this present work.

\VS
In this paper, we extend the FEAST algorithm and theory to tackle non-Hermitian 
generalized eigenproblems.
Similar to the Hermitian case, the non-Hermitian FEAST algorithm takes the
form of standard subspace iteration in
conjunction with the Rayleigh-Ritz procedure (see for example
~\cite{demmel-numerical-linear-algebra}, page 157, or~\cite{saad-eigenvalue-problems-2011}, page
115.)
For non-Hermitian problems, left and right eigenvectors are in general different. 
There are two natural generalizations of subspace iterations to handle this complication. 
A one-sided approach where one focuses
on either the right or left invariant subspace, or a Bauer bi-iteration approach where
both invariant subspaces are targeted simultaneously.
The crucial ingredient is that the subspace iteration here is carried out on an approximate
spectral projector obtained by numerical quadrature. Our analysis shows that the quadrature
approximation perturbs the projector's eigenvalues but not the eigenvectors. 
Consequently, the
convergence of subspace iteration can be established similar to the
approaches shown in~\cite{saad-eigenvalue-problems-2011}, suitably generalized as
the left and right eigenspaces are now different.
By exploring the structure of the generated subspaces,
we show that the Rayleigh-Ritz procedure produces the targeted eigenpairs.
Typical to many large-scale applications,
the target eigenpairs are a small portion of the entire spectrum. In this case,
the dominant work of our algorithm is the quadrature computation which possesses 
multiple levels of parallelism, making 
this an excellent algorithm for high-performance computing.

\VS
This paper aims to show how the various components of non-Hermitian FEAST
fit together, stating the relevant mathematical properties without rigorous
proofs. A detailed numerical analysis similar to~\cite{tang-polizzi-2013-FEAST} 
for the Hermitian case is beyond the scope here. 
In subsequent sections we will describe the numerical-quadrature-based method to compute approximate
spectral projectors, state the convergence properties of subspace iteration 
and the associated Rayleigh-Ritz procedure with this
approximate projector as an accelerator, and present numerical and performance examples.


\section{Overview}
\label{sec:overview}

Throughout this paper, we consider the generalized eigenvalue problem specified
by two $n\times n$ matrices $A$ and $B$, $A \neq 0$ and $B$ invertible. We
assume that $\inverse{B}A$ is diagonalizable with an eigendecomposition
$M \bydef \inverse{B}A = X \, \Lambda \,\inverse{X}$,
$\Lambda$ is a diagonal matrix containing the eigenvalues in some order. $X$ is
a set of corresponding right eigenvectors, $AX = BX\Lambda$. Define
$Y$ by $\ctrans{Y}BX = I$, $\ctrans{Y}$ being the conjugate transpose of $Y$.
$Y$ is a set of corresponding left eigenvectors $\ctrans{Y}A = \Lambda \ctrans{Y}B$
(or $\ctrans{A}Y = \ctrans{B} Y \ctrans{\Lambda}$).
Thus,
\begin{equation}
\inverse{B}A = X\Lambda\ctrans{Y}B, \quad A = BX\Lambda\ctrans{Y}B.
\label{eqn:eigendecomposition}
\end{equation}
It is customary to describe the relationship $\ctrans{Y}BX = I$ as $X$ and $\ctrans{B}Y$
being bi-orthogonal.

\VS
Consider that the eigenvalues of interest, totaling $m$ of them, are 
those that reside inside  a simply connected domain $\domain$ (e.g. disk, ellipse, etc.). 
We further assume that none of the eigenvalues $\lambda\in\eig(\Lambda)$ are on 
the boundary $\boundary$ of $\domain$.
Let $X_\domain$ and $Y_\domain$ be a corresponding set of right and
left eigenvectors, respectively. In particular, $X_\domain$ and $Y_\domain$ are $n \times m$
matrices with $\ctrans{(Y_\domain)} B X_\domain = I_m$. 
Our strategy is motivated by the spectral projectors
onto the invariant subspaces $\spann(X_\domain)$ of $\inverse{B}A$ and
$\spann(Y_\domain)$ of $\invctrans{B}\ctrans{A}$, respectively. In matrix form,
these projectors are $\Rprojector$ and $\Lprojector$.
More specifically, suppose we could compute $\Rprojector\,\vvu$ for any $n$-vector
$\vvu$, then we can apply $\Rprojector$ to a set of random vectors
$U = [\vvu_1,\ldots,\vvu_p]$. Clearly,
$\spann(\Rprojector\,U) \subseteq \spann(X_\domain)$. If it happens that
$\rank(\Rprojector\,U) = \rank(X_\domain)$, then
$\spann(\Rprojector\,U) = \spann(X_\domain)$. 
One can then obtain a basis $\widehat{U}$ for
$\spann(X_\domain)$, for example by performing a rank-revealing factorization. Thus
$\widehat{U}$ must be of the form $\widehat{U} = X_\domain \inverse{W}$ for some
$m \times m$ matrix $W$. Construct a reduced-size eigenproblem $(\widehat{A},\widehat{B})$ where
$
\widehat{A} \bydef \ctrans{\widehat{U}}A\widehat{U}, \quad
\widehat{B} \bydef \ctrans{\widehat{U}}B\widehat{U}.
$
It is easy to see from Equation~(\ref{eqn:eigendecomposition}) that
$\widehat{A} W = \widehat{B} W \Lambda_\domain$. Solving the reduced-size
problem 
for $\Lambda_\domain$ and $W$
yields the eigenvalues of interest and
the eigenvectors $X_\domain$,
which are given by $X_\domain = \widehat{U} W$.

\VS
Similarly, the projector $\Lprojector$ can lead to a basis 
$\widehat{V} = Y_\domain \inverse{Z}$. Construct the reduced-size generalized
eigenproblem $(\widehat{A},\widehat{B})$, 
$
\widehat{A} \bydef \ctrans{\widehat{V}}\ctrans{A}\widehat{V}, \quad
\widehat{B} \bydef \ctrans{\widehat{V}}\ctrans{B}\widehat{V}.
$
Solving $\widehat{A}Z = \widehat{B} Z \Lambda^{H}_\domain$ for
$\Lambda^{H}_\domain$ and $Z$ yields $\Lambda_\domain$ and
$Y_\domain = \widehat{V}\, Z$. 
Finally, if we employ both projectors to obtain basis $\widehat{U}$ and
$\widehat{V}$ for $\spann(X_\domain)$ and $\spann(Y_\domain)$, respectively,
then we can construct
$$
\widehat{A} \bydef \ctrans{V} A U = \invctrans{Z} \Lambda_\domain \inverse{W},
\quad
\widehat{B} \bydef \ctrans{V} B U = \invctrans{Z} \inverse{W}.
$$
The eigenvalues of the reduced problem is $\Lambda_\domain$ and
the right and left eigenvectors are $W$ and $Z$, respectively.

\VS
While the exact spectral projectors are not readily available, we show how
we can approximate them based on rational approximations to a Cauchy integral
via quadrature rules. Applying these approximate projectors is tantamount to solving
multiple independent linear systems each with multiple right-hand-sides --
a procedure that is inherently parallel in multiple levels. Furthermore,
the approximate spectral projectors in fact preserve the invariant subspaces
$\Span(X_\domain)$ and $\Span(Y_\domain)$ exactly.
Consequently, performing 
subspace iteration or Bauer bi-iteration with the approximate spectral
projector becomes numerically effective as well as computationally efficient
in capturing invariant subspaces as well as the associated eigenpairs.
The general flow of the remaining sections is as follows.
In Section~\ref{sec:asp_via_quad}, 
we construct the approximate spectral projectors and analyze their properties.
Section~\ref{sec:subspace_iteration} presents several variants of 
the approximate-spectral-projector-accelerated subspace iteration algorithms
adapted for generalized non-Hermitian eigenvalue problems. 
We state the basic convergence properties of these methods.
We present in Section~\ref{sec:experiments} 
a number of numerical experiments to illustrate the 
theoretical analysis. Scalability results are also presented, supporting our claim that 
this building block is a great addition to the overall toolbox for HPC calculation of non-Hermitian
eigenvalue problems.  In the concluding section, we put our new method in the context 
of other popular existing methods and share our views of future work.


\section{Projection via Quadrature}
\label{sec:asp_via_quad}

We focus first on the right projector $\Rprojector$. Let
$\quadf(\mu) = \sum_{k=1}^q \alpha_k /(\beta_k - \mu)$ be a rational function
in partial fraction form where 
$\beta_k \notin \eig(\Lambda)$ for all the $\beta_k$s.
The standard definition of $\quadf$ applied to the matrix $\surM$ is
\begin{equation}
\quadf(\surM) = \sum_{k=1}^q \alpha_k\,(\beta_k I - \surM)^{-1} = X\,\quadf(\Lambda) \ctrans{Y}B.
\label{eqn:quadf_surM}
\end{equation}
The last equality holds 
because $\surM = X\Lambda\ctrans{Y}B$. Here $\quadf(\Lambda)$ has the obvious meaning of
the diagonal matrix with entries $\quadf(\mu)$. If it happens that
$\quadf(\mu)=1$ for all eigenvalues $\mu \in \domain$ and $\quadf(\mu)=0$
for all eigenvalues $\mu \notin \domain$, then in fact
$\quadf(\surM) = \Rprojector$. In the following, we construct a function $\quadf(\mu)$
such that $\quadf(\mu) \approx 1$ for $\mu \in \domain$ and $\quadf(\mu) \approx 0$ for
$\mu \notin \domain$. Consequently, the resulting $\quadf(\surM)$ approximates the spectral
projector $\Rprojector$.

\VS
Let $\pi(\mu)$ be the complex-valued function defined by the Cauchy integral (in the
counter clockwise direction)
\begin{equation}
\pi(\mu) = \ConInt{ \frac{1}{z - \mu} },
\qquad \mu \notin \boundary.
\label{eqn:cauchy_contour}
\end{equation}
The Cauchy integral theorem shows that $\pi(\mu) = 1$ for $\mu$ inside the 
$\domain$
and $\pi(\mu) = 0$ for $\mu$ outside of $\domain$. 
It is therefore natural to approximate the integral
in Equation~(\ref{eqn:cauchy_contour})
by a quadrature rule. 
To simplify the exposition, we focus on elliptical contours $\boundary$
parameterized by
\begin{equation}
\label{eqn:ellipses}
\phi(t) = c + R\left(
                  \cos\left(\frac{\pi}{2}(1+t)\right) + 
          \eye  a \sin\left(\frac{\pi}{2}(1+t)\right) \right),
\end{equation}
where $c\in\complex$ is the center, with horizontal
and vertical axes of lengths $R> 0$  and $a > 0$, respectively.
We can apply any quadrature rule 
for integrating a function $f(t)$ on $[-1,1]$ 
to obtain an approximation of $\pi(\mu)$. Let
$\int_{-1}^{1} f(t) dt \approx \sum_{k=1}^K w_k f(t_k)$
be a quadrature rule based on
$K$ pairs of $({\rm node}, {\rm weight})$, 
$\{(t_k,w_k) | t_k \in [-1,1], w_k > 0, k = 1, 2, \ldots, K\}$.
\begin{eqnarray}
\pi(\mu) & = & \ConInt{ \frac{1}{z-\mu} }, \nonumber \\
 & = & \frac{1}{2\pi\eye} \int_{-1}^3 \frac{\phi'(t)}{\phi(t) - \mu} dt, \nonumber \\
 & = & \frac{1}{2\pi\eye}
       \int_{-1}^1 \left(
       \frac{\phi'(t)}{\phi(t)-\mu} + \frac{\phi'(2-t)}{\phi(2-t)-\mu}
       \right).
\label{eqn:pi}
\end{eqnarray}
Applying the quadrature rule
\footnote{
See~\cite{sakurai-sugiura-2003} for a different
application of numerical quadrature to eigenvalue problems.}
yields $\quadf(\mu)\approx\pi(\mu)$,
\begin{eqnarray}
\quadf(\mu)
 & \bydef &
\sum_{k=1}^K \frac{w_k}{2\pi\eye} 
\left( \frac{\phi'(t_k)}{\phi(t_k)-\mu} +
       \frac{\phi'(\tilde{t}_k)}{\phi(\tilde{t}_k)-\mu} 
\right), \quad \tilde{t}_k = 2-t_k, \nonumber \\
 & = & \sum_{k=1}^q \sigma_k (\phi_k - \mu)^{-1}, \quad
       \sigma_k, \phi_k \in \complex.
\label{eqn:quadrature_for_pi}
\end{eqnarray}
$\quadf(\surM)$ acting on any $n\times p$ matrix $U$, 
$\quadf(\surM)\,U$ is
\begin{equation}
 \sum_{k=1}^q \sigma_k (\phi_k I - \surM)^{-1} \, U
  = 
 \sum_{k=1}^q \sigma_k (\phi_k B - A)^{-1} \, (BU).
\label{eqn:Rasp}
\end{equation}
The operation involves solving $q$ linear systems each with $BU$ as the
right-hand-side. If the $K$-point quadrature rule is such that
neither $-1$ nor $1$ is used as nodes, then $q = 2K$. If $t_1 = -1$
and $t_K = 1$, then $q = 2(K-1)$.
Solutions of multiple independent linear systems for multiple right hand sides
make $\quadf(\surM)\,U$ a kernel operation with rich parallelism. 
Furthermore, $\quadf(\surM)\,U$ is numerically effective, as we now explain.

\VS
We examine the ratio of
$|\quadf(\mu)/\quadf(\mu')|$ for $\mu'\in\domain$ and $\mu\notin\domain$.
To this end, if suffices to study the reference $\quadf$ function
$\quadf_\reference$ for the domain $\domain$ that centers at the origin, with
$R = 1$ because the $\quadf(\mu)$ function for an ellipse of a same ``$a$''
parameter but centered at $c$ with ``radius'' $R'$ is simply given
$\quadf(\mu)  = \quadf_\reference((\mu-c)/R')$.

\begin{figure}
\centering 
\includegraphics[width=3.375in]{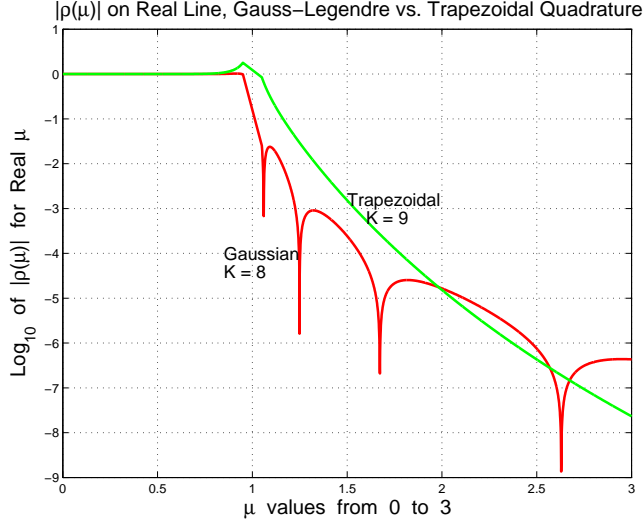}
\caption{For Hermitian eigenproblems, one only needs a $\quadf_\reference(\mu)$
function that behaves well on the real line, namely, is close to 1 inside $[-1,1]$ and 
small outside. Because of symmetry, this figure shows $|\quadf_\reference(\mu)|$ only
for $\mu \ge 0$. For comparable computational cost of applying $\quadf(\surM)$ to vectors
$U$, Gauss-Legendre performs somewhat better than trapezoidal rule does.} 
\label{figure:rho_on_real}
\end{figure}

\VS
To underline the difference between Hermitian and non-Hermitian problems,
an effective quadrature rule for the former requires 
$|\quadf_\reference(\mu)| \approx 1$ for $\mu \in \domain$ and
$|\quadf_\reference(\mu)| \ll 1$ for $\mu \notin \domain$ 
\emph{only} for $\mu$ on the real line. Figure~\ref{figure:rho_on_real}
shows $\log_{10}|\quadf_\reference(\mu)|$ for real-valued $\mu$
for a Gauss-Legendre (with $K = 8$) and a trapezoidal rule (with $K=9$).
The precipitous drop of $|\quadf_\reference(\mu)|$ for $\mu$ outside
of $[-1,1]$ signifies the effectiveness of quadrature-based approximate
spectral projections.

\begin{figure}
\centering 
\includegraphics[width=3.375in]{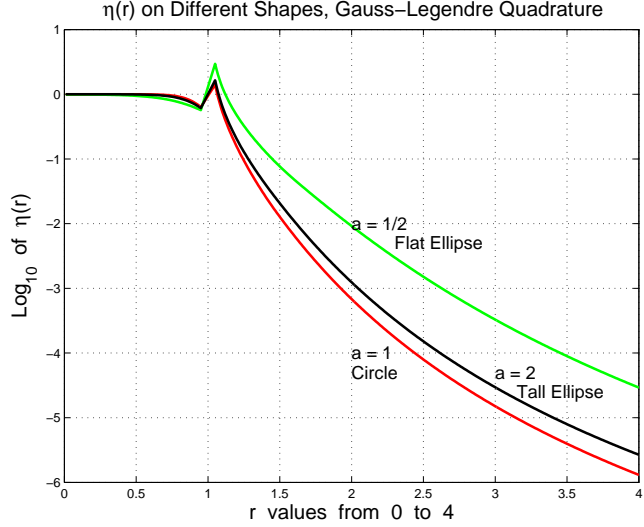}
\caption{
These are the $\AI(r)$ functions corresponding to Gauss-Legendre
quadrature with $K=8$ nodes on $[-1,1]$. We exhibit the result
for three different elliptical domains. For simplicity, we employ circular
domains for the rest of the paper, but different types of domains 
can be used. See further discussions in Section~\ref{sec:conclusion}.}
\label{figure:eta_on_different_shapes}
\end{figure}

\begin{figure}
\centering 
\includegraphics[width=3.375in]{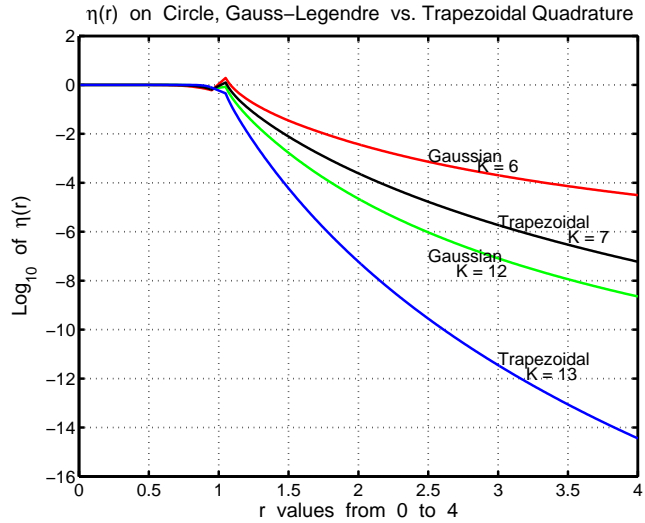}
\caption{
This figure compares Gauss-Legendre quadrature to trapezoidal
rule on a circular domain. Because trapezoidal rule uses both
$-1$ and $1$ as nodes on the integration interval $[-1,1]$ while
Gauss-Legendre uses neither, a $K$-node Gauss-Legendre and
a $K+1$-node trapezoidal both require solving $2K$ linear systems
when applying the spectral projector $\quadf(\surM)$ to vectors $U$.
The figure suggests that trapezoidal rule works better in general
for non-Hermitian problems.}
\label{figure:eta_circle_GL_vs_TR}
\end{figure}

\VS
For non-Hermitian problems, $\quadf_\reference(\mu)$ has to ``behave well''
for $\mu$ in the complex plane. Consequently, for a given quadrature rule,
we evaluate $\quadf_\reference(\mu)$ at level curves similar to the 
boundary $\boundary$:
$$
\mu(r,t) = r \left[        \cos\left(\frac{\pi}{2}(1+t)\right) + 
                        \eye a \sin\left(\frac{\pi}{2}(1+t)\right)
                 \right].
$$
At each $r$ below 1, $0 \le r \le 1-\delta$ ($\delta$ set to $0.01$), we record
the minimum of $|\quadf_\reference|$ over the level curve, and
at each $r \ge 1+\delta$, we record the maximum. That is, we examine
the function
$$
\AI(r) \bydef \left\{
\begin{array}{c l}
\min_t |\quadf_\reference(\mu(r,t))| & \hbox{for $0 \le r \le 1-\delta$,} \\
\max_t |\quadf_\reference(\mu(r,t))| & \hbox{for $1 + \delta \le r$.}
\end{array}
\right.
$$
The function $\AI(r)$ serves as an indicator. An $\AI(r)$ that is close to 1 for $r < 1$
and very small for $r > 1$ corresponds to an approximate spectral projector
that preserves the desired eigenspace well while attenuating the unwanted eigencomponents
severely. Figure~\ref{figure:eta_on_different_shapes} shows three $\AI(r)$ functions, 
in logarithmic scale, corresponding to Gauss-Legendre quadrature ($K=8$) 
on three different shapes of ellipses. 
Figure~\ref{figure:eta_circle_GL_vs_TR} shows different $\AI(r)$ functions,
in logarithmic scale, corresponding to Gauss-Legendre and trapezoidal rules at
different choices of $K$. The domain is set to be a circle. It is interesting to note that
while Gauss-Legendre is in general a better choice for Hermitian problems 
(as Figure~\ref{figure:rho_on_real} suggests), 
trapezoidal rule seems to fare  better for
non-Hermitian problems.\footnote{
Assuming no information of the eigenvalues' distribution is available
a priori.}

\VS
We note that the left projector can be approximated in a similar fashion.
In particular,
\begin{equation}
\invctrans{B}\ctrans{[\quadf(\surM)]}\ctrans{B} =
Y\ctrans{[\quadf(\Lambda)]}\ctrans{X}\ctrans{B} \approx \Lprojector.
\label{eqn:quadf_surMconj}
\end{equation}
Applying this approximate left projector on a $n\times p$ matrix $V$ is
\begin{equation}
\invctrans{B}\ctrans{[\quadf(\surM)]}\ctrans{B} \, V
  = 
 \invctrans{\left[\sum_{k=1}^q \sigma_k (\phi_k B - A)\right]} \, (\ctrans{B}\,V).
\label{eqn:Lasp}
\end{equation}
This involves solving $q$ conjugate-transposed system $\ctrans{[\sigma_k(\phi_k B - A)]}$
each with $\ctrans{B}V$ as right-hand-side.

\section{Non-Hermitian FEAST}
\label{sec:subspace_iteration}

Equation~(\ref{eqn:quadf_surM}) shows that 
$(\quadf(\lambda_j),\vvx_j)$ is a right eigenpair of
$\quadf(\surM) \vvx_j = \quadf(\lambda_j) B \vvx_j$
for every right eigenpair $(\lambda_j,\vvx_j)$ of
$A\vvx_j = \lambda_j B \vvx_j$. Moreover $\eig(\quadf(\Lambda_\domain))$ are
among the most dominant eigenvalues $\eig(\quadf(\Lambda))$ of the approximate
right projector $\quadf(\surM)$. It is easy to see that corresponding
properties hold for the approximate left projector.

\VS
Since $\eig(\quadf(\Lambda_\domain))$ are the dominant eigenvalues
of $\quadf(\surM)$, subspace iteration with $\quadf(\surM)$ is effective in capturing
the invariant subspace $\Span(X_\domain)$. A Rayleigh-Ritz projection of the original
eigenproblem would then allow us to obtain $\Lambda_\domain$
and $X_\domain$. This use of numerical-quadrature-based
approximate spectral projector to accelerate subspace iteration followed by Rayleigh-Ritz
is the essence of the Hermitian FEAST 
algorithm~\cite{polizzi-2009,tang-polizzi-2013-FEAST}.
We offer here two generalizations to non-Hermitian problems. The first,
Algorithm R-FEAST, uses only one projector\footnote{
We use the right projector here, but it is obvious how a left-projector
variant would work.}; the second, Algorithm Bi-FEAST, uses both.

\begin{algorithm}[h!]
{\it Algorithm R-FEAST} (Right-Projector FEAST)
\begin{algorithmic}[1]
\State Pick random $\iterate{U}{0} \in \complex^{n\times p}$. Set $k \gets 1$.
\Repeat
\State $\widehat{U} \gets (\quadf(\surM) \cdot \iterate{U}{k-1})$
\State $\widehat{A} \gets \ctrans{\widehat{U}} A \widehat{U}, \quad
        \widehat{B} \gets \ctrans{\widehat{U}} B \widehat{U}$. 
\State Solve $\widehat{A} W = \widehat{B} W \widehat{\Lambda}_{(k)}$ for
       $\widehat{\Lambda}_{(k)}$ and $W$.
\State $\iterate{U}{k} \gets \widehat{U}\cdot W$.
\State $k \gets k + 1$
\Until {Appropriate stopping criteria}
\end{algorithmic}
\end{algorithm}

\begin{algorithm}[h!]
{\it Algorithm Bi-FEAST} (Bi-iteration FEAST)
\begin{algorithmic}[1]
\State Pick random $\iterate{U}{0} \in \complex^{n\times p}$. 
\State Get $\iterate{V}{0}$ where $\ctrans{(\iterate{V}{0})}B\iterate{U}{0} = I_p$. Set $k \gets 1$.
\Repeat
\State $\widehat{U} \gets (\quadf(\surM)\cdot \iterate{U}{k-1}$
\State $\widehat{V} \gets \invctrans{B}\ctrans{\quadf}(\surM)\ctrans{B}\cdot \iterate{V}{k-1}$.
\State $\widehat{A} \gets \ctrans{\widehat{V}} A \widehat{U}, \quad
        \widehat{B} \gets \ctrans{\widehat{V}} B \widehat{U}$. 
\State Solve $\widehat{A}  = \widehat{B}\, W \, \widehat{\Lambda}_{(k)}\, \ctrans{Z}\widehat{B}
       \; \hbox{for
       $\widehat{\Lambda}_{(k)},\; W,\; {\rm and}\; Z$}$.
\State $\iterate{U}{k} \gets \widehat{U}\cdot W$, $\iterate{V}{k} \gets \widehat{V}\cdot Z$.
\State $k \gets k + 1$
\Until {Appropriate stopping criteria}
\end{algorithmic}
\end{algorithm}

\noindent
We state without proofs several key properties of these two algorithms that correspond to 
generalizations of corresponding theorems in~\cite{tang-polizzi-2013-FEAST}.
Number the $\gamma_j$s (the eigenvalues of $\quadf(\surM)$) so that
$$
|\gamma_1| \ge |\gamma_2| \ge \cdots \ge |\gamma_n|,
$$
and number the eigencomponents of $(A,B)$ accordingly:
$$
\gamma_j = \quadf(\lambda_j), \;
A \vvx_j = \lambda_j B \vvx_j, \; \ctrans{\vvy}_j A = \lambda_j \ctrans{\vvy}_j B, 
\quad j = 1, 2, \ldots, n.
$$
We use the following notations: For integer $\ell$,
$1 \le \ell \le n$,
\begin{equation*}
\begin{array}{l  l}
\XTop{\ell}  =    [\vvx_1, \vvx_2, \ldots, \vvx_\ell],                         &
\XBot{\ell}  =    [\vvx_{\ell+1}, \vvx_{\ell+2}, \ldots, \vvx_n],              \\
\YTop{\ell}  =    [\vvy_1, \vvy_2, \ldots, \vvy_\ell],                         &
\YBot{\ell}  =    [\vvy_{\ell+1}, \vvy_{\ell+2}, \ldots, \vvy_n],              \\
\LambdaTop{\ell}  =   \diag(\lambda_1, \lambda_2, \ldots, \lambda_\ell),  &
\LambdaBot{\ell}  =   \diag(\lambda_{\ell+1}, \lambda_{\ell+2}, \ldots, \lambda_n), \\
\GammaTop{\ell}  =   \diag(\gamma_1, \gamma_2, \ldots, \gamma_\ell),  &
\GammaBot{\ell}  =   \diag(\gamma_{\ell+1}, \gamma_{\ell+2}, \ldots, \gamma_n). 
\end{array}
\end{equation*}

\VS
Figures~\ref{figure:eta_on_different_shapes} and~\ref{figure:eta_circle_GL_vs_TR} are illustrative
of the general properties of quadrature-based functions $\quadf(\mu)$. 
In general $|\quadf(\mu)| \approx 1$ for the $m$ eigenvalues $\lambda \in \eig(\Lambda_\domain)$ 
of interest. They are among the dominant eigenvalues which include those, if exist, outside
of $\domain$ but close to the boundary $\boundary$. Thus, there is an integer
$m' \ge m$, $m' \approx m$, such that
$\eig(\Lambda_\domain) \subseteq \{\lambda_1,\lambda_2,\ldots,\lambda_{m'}\}$ and
$|\quadf(\lambda_j)/\quadf(\lambda_i)| < 1$ for all $i\le m'$, $j > m'$.
Below are the relevant properties of Algorithms R-FEAST and Bi-FEAST under the assumptions
$p \ge m'$, $|\gamma_{p+1}/\gamma_{m'}| \ll 1$ and other moderate technicalities.

\begin{enumerate}
\item
There is a constant $\alpha$ such that for each iteration $k$ and
$j=1,2,\ldots,p$, elements of the form
$\vvx_j + \XBot{p} \vve_j^{(k)}$ and of the form
$\vvy_j + \YBot{p} \vvf_j^{(k)}$ exist in $\Span(\iterate{U}{k})$
and  $\Span(\iterate{V}{k})$, respectively, such that
$\Enorm{\vve_j^{(k)}} ,\Enorm{\vvf_j^{(k)}} \le \alpha \eps^k$
where $\eps = |\gamma_{p+1}/\gamma_{m'}|$.
In other words, as long as $p$ is chosen big enough so that $|\gamma_{p+1}|$
is small (see Figure~\ref{figure:eta_on_different_shapes} for example),
as iterations proceed,
there are elements in $\Span(\iterate{U}{k})$ that are close to
$\Span(X_\domain)$, and elements in $\Span(\iterate{V}{k})$ that are close
to $\Span(Y_\domain)$.

\item
For Algorithm R-FEAST, as iterations proceed, the eigenvalues of 
$\inverse{\widehat{B}}\widehat{A}$ (see Step 4) are the same as those of
a matrix of the form
$$
\left[
\begin{array}{c | c}
A'_{11} & A'_{12} \\  \hline
A'_{21} & A'_{22}
\end{array}
\right],
\quad
\hbox{$A'_{11}$ is $m'\times m'$,}
$$
where $\Enorm{A'_{11} - \TSec{\Lambda}{m'}} = O(\eps^k)$,
and $\Enorm{A'_{21}} = O(\eps^k)$, $\eps = |\gamma_{p+1}/\gamma_{m'}|$.
In particular, there are $m'$ eigenvalues $\hat{\lambda}_j \in \eig(\iterate{\widehat{\Lambda}}{k})$,
$j=1,2,\ldots,m'$,
and the corresponding column vectors $\vvu_j^{(k)}$ of $\iterate{\widehat{U}}{k}$ that satisfy
$$
|\hat{\lambda}_j - \lambda_j| = O(\eps^k),
\quad
\Enorm{(A - \widehat{\lambda}_j B) \vvu_j^{(k)}} = O(\eps^k).
$$

\item
For Algorithm Bi-FEAST, as iterations proceed, the eigenvalues of 
$\inverse{\widehat{B}}\widehat{A}$ (see Step 4) are the same as those of
a matrix of the form
$$
\left[
\begin{array}{c | c}
A'_{11} & A'_{12} \\  \hline
A'_{21} & A'_{22}
\end{array}
\right],
\quad
\hbox{$A'_{11}$ is $m'\times m'$,}
$$
where $\Enorm{A'_{11} - \TSec{\Lambda}{m'}} = O(\eps^{2k})$,
and $\Enorm{A'_{21}}, \Enorm{A'_{12}} = O(\eps^k)$, $\eps = |\gamma_{p+1}/\gamma_{m'}|$.
In particular, there are $m'$ eigenvalues $\hat{\lambda}_j \in \eig(\iterate{\widehat{\Lambda}}{k})$,
$j=1,2,\ldots,m'$, 
and the corresponding column vectors $\vvu_j^{(k)}$ of $\iterate{\widehat{U}}{k}$,
$\vvv_j^{(k)}$ of $\iterate{\widehat{V}}{k}$
that satisfy
$
|\hat{\lambda}_j - \lambda_j| = O(\eps^{2k})
$,
$$
\Enorm{(A - \widehat{\lambda}_j B) \vvu_j^{(k)}},
\quad \Enorm{\ctrans{(A  - \widehat{\lambda}_j B)} \vvv_j^{(k)}}
 = O(\eps^k).
$$
In the absence of ill-conditioning, the discussions of which is omitted here, Bi-FEAST offers
faster convergence of eigenvalues (but not the residuals) compared with R-FEAST, albeit at a higher
computation cost per iteration -- needing to solve the conjugated systems as well. We must
mention that R-FEAST is inherently more stable, especially if orthogonalization is applied to
$\widehat{U}$ between Steps 4 and 5.
\end{enumerate}


\section{Numerical Experiments}
\label{sec:experiments}

\begin{figure}[h!]
\centering 
\includegraphics[width=3.375in]{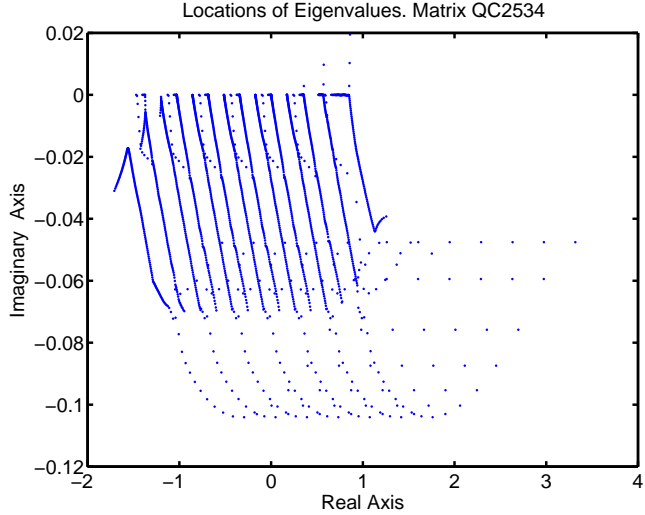}
\caption{Eigenvalues of the matrix QC2534 from the NEP collection.
Eigenvalues distribution of QC234 bears a resemblance.}
\label{figure:QC_eig}
\end{figure}

Experiments given in Sections~\ref{example:R-FEAST} through~\ref{example:domain_shapes}
illustrate the various properties of R-FEAST and Bi-FEAST. They use the matrices
QC324 and QC2534 from the NEP 
collection in~\cite{bai-et-al-1996}.
These are standard non-Hermitian eigenvalue problems ($B=I$) that arise in
quantum chemistry~\cite{chu-qc2534-1991}. 
These two matrices are similar in properties but differ in size.
Figure~\ref{figure:QC_eig} profiles the location of the eigenvalues in the complex plane.
Experiments on QC324 and QC2534 are run in matlab.
It turns out that the $m$ eigenvalues in the domain correspond to the most dominant
eigenvalues in the approximate projectors, thus $m=m'$. Subspace dimensions $p$ in
these experiments are chosen moderately bigger than $m$. 
Section~\ref{sec:conclusion} will discuss how $p$ is set in practice.
The remaining experiments show FEAST applied to actual scientific applications,
run on different computer clusters.

\VS
During the iterations of FEAST, we monitor the $p$ 
eigenpairs computed from the reduced system
(in Step 5 of R-FEAST, for example). A particular (right) eigenpair 
$(\widehat{\lambda}_j,\vvu_j)$ is considered a candidate
if $\widehat{\lambda}_j \in \domain$ and 
$
\Enorm{A \vvu_j - \widehat{\lambda}_j B  \vvu_j} / \Enorm{\vvu_j}
$
is reasonably small, typically, $\le 10^{-4}$. 
Specifically, we track convergence of
$$
{\rm Res}_{(k)} \bydef
\max
\Enorm{A \vvu_j - \widehat{\lambda}_j B  \vvu_j} / \Enorm{\vvu_j}, \;\;
{\rm trace}_{(k)} \bydef
\sum \widehat{\lambda}_j,
$$
over the candidates $\widehat{\lambda}_j$s.

\subsection{Simple Convergence of R-FEAST}
\label{example:R-FEAST}

We illustrate the most basic convergence properties with the small (dimension 324)
matrix QC324.  The domain $\domain$ chosen is the disk of radius $0.01$ centered on the 
real axis at $-0.5$, containing 
$m = 8$ eigenvalues. We employ Gauss-Legendre quadrature and picked $p=8$. 
With this choice, $\log_{10}|\quadf(\lambda_{p+1})/\quadf(\lambda_m)| = -0.93$.  
The table here exhibits the expected behavior from both R-FEAST and Bi-FEAST. 
The eigenvalue and residual convergence rate are linear at roughly $0.9$ digits
per iteration, except that eigenvalues in Bi-FEAST converge as fast as $2\times 0.9$
digits per iteration.

\VS
\begin{center}\begin{tabular}{c | r r | r r}
\multicolumn{5}{c}{
$p=m=8$, Gauss-Legendre with $K=8$} \\

 & \multicolumn{2}{c|}{$\log_{10}|\hbox{change in trace}|$}
 & \multicolumn{2}{c}{$\log_{10}(\hbox{max of residual})$} 
 \\
Iter.  & \RFEAST & \BFEAST &  \RFEAST  & \BFEAST \\ \hline
 4  & -5.4\pk & -0.0\pk & -5.3\pk & -4.8\pk    \\
 5  & -6.4\pk & -8.6\pk & -6.2\pk & -5.7\pk    \\
 6  & -7.3\pk & -10.5\pk & -7.1\pk & -6.6\pk    \\
 7  & -8.2\pk & -12.4\pk & -8.1\pk & -7.6\pk    \\
 8  & -9.2\pk & -14.3\pk & -9.0\pk & -8.5\pk    \\
 9  & -10.2\pk & -14.4\pk & -9.9\pk & -9.4\pk    \\
10  & -11.1\pk & -14.5\pk & -10.9\pk & -10.3\pk    \\
11  & -12.1\pk & -15.1\pk & -11.8\pk & -11.3\pk    \\
12  & -13.1\pk & -14.8\pk & -12.7\pk & -12.2\pk    \\
13  & -14.1\pk & -14.8\pk & -13.7\pk & -13.1\pk    \\
14  & -14.7\pk & -14.5\pk & -14.6\pk & -14.1\pk 
\end{tabular}
\end{center}

\subsection{R-FEAST and Bi-FEAST}
\label{example:R-FEAST-Bi-FEAST}
 
This example illustrates the sensitive nature of Bi-FEAST. We have seen in the
previous example that Bi-FEAST can offer a faster convergence on
the eigenvalues. But as discussed in Section~\ref{sec:subspace_iteration}, 
Bi-FEAST is more
sensitive to the conditioning of the eigenvalues. This is the case
for the matrix QC2534 when the region is chosen to be 
the disk of radius $0.01$ centered on the real axis at $0.85$, containing 10 eigenvalues.
With $p$ set to $p = m+5 = 15$, Gauss-Legendre quadrature with $K=8$ 
yields $\log_{10}|\quadf(\lambda_{p+1})/\quadf(\lambda_m)| = -2.61$.  
The condition of the eigenvalues, however, are poor: the
products $\Enorm{\vvx_j}\, \Enorm{\vvy_j}$ are of the order of $10^{11}$.
The table here shows that indeed the eigenvalues cannot be resolved to be much better than 5 or 6 digits.
R-FEAST is able to deliver small residuals, while Bi-FEAST is hampered by the poor conditioning, as it
is difficult to maintain bi-orthogonality between the $\vvx_j$ and $\vvy_j$ to 
full machine precision, precisely because their norms are large.

\VS
\begin{center}\begin{tabular}{c | r r | r r}
\multicolumn{5}{c}{
$p=m+5=15$, Gauss-Legendre with $K=6$} \\

 & \multicolumn{2}{c|}{$\log_{10}|\hbox{change in trace}|$}
 & \multicolumn{2}{c}{$\log_{10}(\hbox{max of residual})$} 
 \\
Iter.  & \RFEAST & \BFEAST &  \RFEAST  & \BFEAST \\ \hline
 2  & 0.6\pk & 1.0\pk & -8.3\pk & -4.1\pk    \\
 3  & -0.0\pk & -5.3\pk & -11.4\pk & -5.8\pk    \\
 4  & -5.2\pk & -5.3\pk & -14.0\pk & -5.9\pk    \\
 5  & -6.8\pk & -5.4\pk & -14.2\pk & -6.2\pk    \\
 6  & -6.9\pk & -5.3\pk & -14.2\pk & -5.9\pk    \\
 7  & -7.6\pk & -5.3\pk & -14.2\pk & -6.1\pk    \\
 8  & -6.9\pk & -5.6\pk & -14.2\pk & -6.0\pk    \\
 9  & -6.6\pk & -5.4\pk & -14.3\pk & -6.0\pk    \\
10  & -6.8\pk & -5.7\pk & -14.1\pk & -5.8\pk 
\end{tabular}
\end{center}

\subsection{Different Quadratures}
\label{example:quadratures}

Figure~\ref{figure:eta_circle_GL_vs_TR} in Section~\ref{sec:asp_via_quad} suggests that
trapezoidal rule may work better in general. 
This example is consistent with this view, but illustrates some subtlety. 
Figure~\ref{figure:eta_circle_GL_vs_TR} depicts minimal convergence rate. Depending on the exact
location of the eigenvalues, which is problem specific, a quadrature with a lower minimal convergence
rate may actually still converge faster.
Here we compute the eigenvalues of QC2534 that reside inside the disk of radius $0.02$, centered
on the real line at $-0.17$, which contains 28 eigenvalues. At each of two different settings,
the table below exhibits the residual convergence for both Gauss-Legendre and trapezoidal quadrature.
The behavior below is consistent with the actual values of $|\quadf(\lambda_{p+1})/\quadf(\lambda_m)|$.

\VS
\begin{center}\begin{tabular}{c | r r | r r}
\multicolumn{5}{c}{
R-FEAST, 
Gauss-Legendre(GL) vs. Trapezoidal(TR)}
\\
 &
\multicolumn{4}{c}{
$\log_{10}(\hbox{max of residual})$
}
\\
 & \multicolumn{2}{c|}{$p=m+3=31$}
 & \multicolumn{2}{c}{$p=m+6=34$}
 \\
Iter.  & {\scriptsize GL-8 nodes} & {\scriptsize TR-9 nodes} 
       & {\scriptsize GL-8 nodes} & {\scriptsize TR-9 nodes} \\ \hline
 2  & -4.1\pk & -4.0\pk & -4.6\pk & -5.7\pk    \\
 3  & -5.6\pk & -5.4\pk & -6.4\pk & -8.2\pk    \\
 4  & -7.0\pk & -6.8\pk & -8.3\pk & -11.3\pk    \\
 5  & -8.7\pk & -8.2\pk & -10.1\pk & -13.8\pk    \\
 6  & -10.9\pk & -9.4\pk & -11.9\pk & -14.2\pk    \\
 7  & -12.9\pk & -10.6\pk & -13.7\pk & -14.3\pk    \\
 8  & -14.1\pk & -11.9\pk & -14.4\pk & -14.4\pk 
\end{tabular}\end{center}

\begin{figure}[h!]
\centering 
\includegraphics[width=3.375in,height=2.75in]{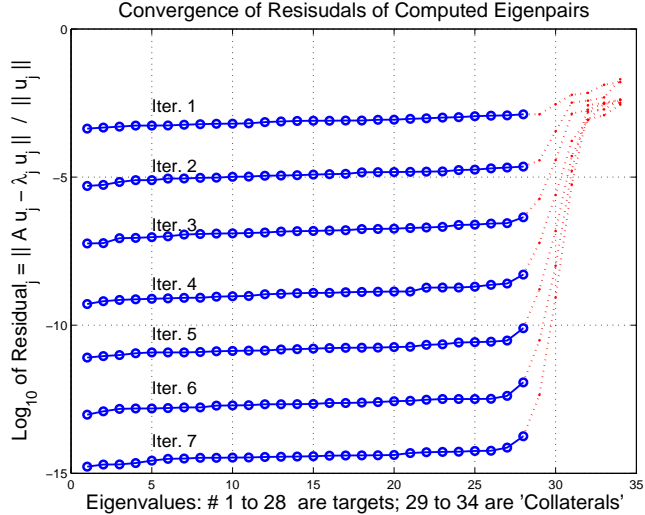}
\caption{Convergence of residuals: targeted and the ``collaterals''.
Residuals are sorted to give a more tidy picture.}
\label{figure:residuals_with_collaterals}
\end{figure}

\VS
The typical convergence pattern of the residuals is as follows. 
The subspace dimension $p$ is in general
bigger than the number of eigenvalues inside the targeted domain. Some of the residuals
that are not targeted (we usually call them collaterals) will converge slowly,
or not at all. Figure~\ref{figure:residuals_with_collaterals} displays the residuals
of our current QC2534 test using Gauss-Legendre with $p$ set to $m+6$. 
Notice that the 28 targeted residuals
converge linearly at the expected rate. Convergence of the collaterals are much
slower, and some not at all.

\subsection{General Domain Shapes}
\label{example:domain_shapes}
\begin{figure}[h!]
\centering 
\includegraphics[width=3.375in]{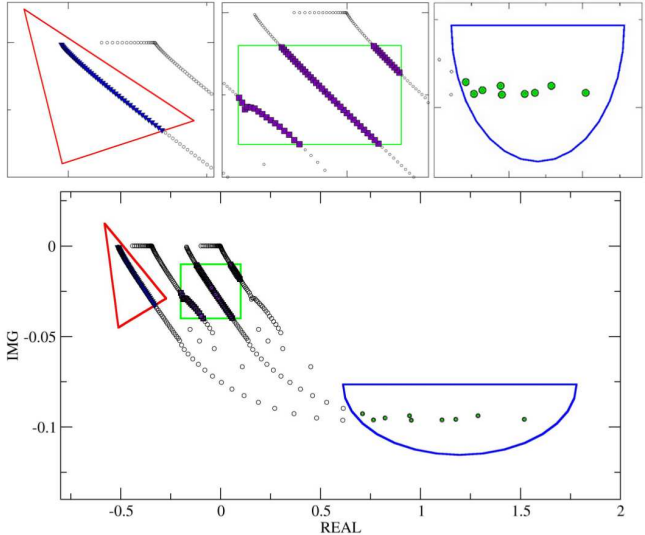}
\caption{Bi-FEAST is run on three different regions to illustrate that
general shapes can be supported easily. The triangle, square, and
a semicircle regions are bounded by 3, 4, and 2 segments of curves (or lines)
parameterized separately and ``glued'' together.}
\label{figure:custom-contours}
\end{figure}

FEAST can obviously be applied on general domain shapes as one only need to change
the parameterization function $\phi(t)$ in Equation~\ref{eqn:ellipses} appropriately.
Using QC2534 still, we compute the spectrum in three different regions, as shown in
Figure~\ref{figure:custom-contours}, using the simple-to-use trapezoidal rule. The table
below summarizes the results of running Bi-FEAST.

\VS
\begin{center}\begin{tabular}{l r | c c c }
 &  & Triangle  &   Square    &   Semicircle \\ \hline

\# eigenvalues &  $m$       &   45   &  64  &  9   \\
subspace dim.  &  $p$       &   80   &  100 &  80  \\
\# quadrature nodes & $K$   &   24   &  32  &  16  \\
\multicolumn{5}{c}{convergence, measured in digits per iteration} \\
\multicolumn{2}{l|}{eigenvalues} &  4  &  2  &  5 \\
\multicolumn{2}{l|}{residuals}   &  2  &  1  &  2.5 \\
\end{tabular}\end{center}

\subsection{An Electronic Structure Problem}
\label{example:application1}

We study the structure of Benzene molecule via a  
finite element (FEM) discretization of the Kohn-Sham equation.
An all-electron potential framework described in~\cite{levin-zhang-polizzi-2012}
is used here. 
Complex potential interfaces at the edges of the 
computational domain~\cite{lehtovaara-havi-puska-2011} are then
added, resulting in complex symmetric eigenvalue problems.
We exhibit here two test cases of quadratic and cubic finite elements, denoted
as FEM-Q and FEM-C.
The eigenproblems are complex symmetric and in generalized form
$Ax = \lambda Bx$.
Algorithm Bi-FEAST is run with Gauss-Legendre, $K=12$, on the
Phoenix cluster at University of Massachusetts
consisting of multiple nodes of \Intel \Xeon X5550 2.66GHz processor, 8 cores per node. 
We use only one target circular domain of radius $3.364 \times 10^{-18}$
centered on the real line at $-2.403\times 10^{-17}$. 
Parallelism on the MPI-process level is exploited by the 24 ($K=12$)
linear systems to be solved.
When there are enough MPI processes, each linear system can be factored just once for
the entire iterative process. This benefit can be seen from the 6-core result of FEM-C
below.  Within each MPI process, parallelism is utilized by
the direct sparse solver Pardiso from \Intel Math Kernel Library version 10.36.
With subspace dimension set to $p=300$, Bi-FEAST
converges on the eigenvalues within 4 and 5 iterations, on FEM-Q and FEM-C,
respectively. The table presented here emphasizes the ability of Bi-FEAST to utilize
the multiple nodes of the system, expressed in ``Efficiency.''

{\scriptsize
\begin{center}\begin{tabular}{l || c c c c  |  c c c c }
&
\multicolumn{4}{c|}{
FEM-Q $n = 149,192$
} &
\multicolumn{4}{c}{
FEM-C $n = 1,165,485$ 
} \\
&
\multicolumn{4}{c|}{
$m=149$, $p=300$
} &
\multicolumn{4}{c}{
$m=181$, $p=300$
} \\ \hline
Time    &  & & & & & & & \\
(Ksec)      &  0.27 &  0.15   &   0.10  &  0.05 
            &  1.65  & 0.86   &   0.61  &  0.25   \\
  &  & & & & & & & \\ 
Eff.     &  & & & & & & & \\ 
(\%)
           &  100   &  92   &   88   &  85  
           &  100   &  96   &   91   &  110   \\  \hline
& 1 & 2 & 3 & 6 & 1 & 2 & 3 & 6  \\
& \multicolumn{8}{c}{Number of Nodes, 2 MPI Process/Node}   
\end{tabular}\end{center}
}

\subsection{Quantum Transport}
\label{example:application2}

We now seek to obtain the quantum bound-states of a Benzene molecule sandwiched between two electrodes.
One can show that an exact derivation of the boundary conditions of the system can give rise to a 
Hermitian but quadratic (non-linear) eigenvalue problem~\cite{shao-porod-lent-1995}.
From this model, however, one can formulate a more practical  
linear companion problem but twice larger and non-Hermitian.
The determination of the resonant states, that is, solution of this non-Hermitian problem is 
essential in quantum transport theory~\cite{polizzi-abdallah-vanbesien-lippens-2000}. 
The eigenvalues of interest are located close to the real axis.

\VS
A non-Hermitian problem of size $n=98,384$ is solved on a cluster named Endeavor, which is part of
the computing infrastructure of Intel Corporation. The search domain is a circle 
in $\vvk$-space chosen to contain the energy resonances
of the non-linear problem. These resonances correspond to eigenvalues with
tiny imaginary parts, in the range $-0.99$eV to $2$eV.
Subpsace dimension set to 150, number
of quadrature nodes is set to 24, resulting in 48 linear systems in the evaluation
of the approximate spectral projector. However, only 24 matrix factorizations are needed
because search domain boundary is symmetric with respect to the real axis. In this
setting, eigenvalues converge at the second iteration using Bi-FEAST. The following table
shows the compute time when Bi-FEAST is run on multiple compute nodes, each node using 16 cores
of third generation Intel\textsuperscript{\textregistered} Xeon\textsuperscript{\textregistered} processor.

\VS

\VS
{\scriptsize
\begin{center}
\begin{tabular}{l | c c c c c c c c }

                & \multicolumn{8}{c}{Number of compute nodes, 16 cores/node} \\
                &  1  &  2  &  3  &  4  &  6  &  8  &  12  &   24   \\ \hline
Total time      &  \\
in seconds      &  366  &  154  &  106  &  77.0  &  54.4 &  44.1  &  34.8  & 24.6  \\
                &  \\
Efficiency (\%) &  100  &  118  &  115  &  118   &  112  & 104    &  90    & 62  

\end{tabular}\end{center}
}


\section{Conclusion}
\label{sec:conclusion}

In the paper, we have introduced a new non-Hermitian eigensolver with rich inherent
parallelism.
This paper establishes the theory behind generalizing 
the FEAST solver for Hermitian problems~\cite{polizzi-2009,tang-polizzi-2013-FEAST} 
to non-Hermitian problems in two flavors. 
Bi-FEAST is the ``bullish-but-riskier" sibling of the 
more conservative R-FEAST. For well-conditioned problems, 
Bi-FEAST offers faster convergence of eigenvalues;
R-FEAST, however, is just as fast in producing small residuals. 
Both are useful and complement each other.
We note here that Bi-FEAST was experimented in~\cite{laux-2012} by Laux, but
without theoretical explanation.

\VS
FEAST has a number of signature features. By nature it works equally well
regardless whether the targeted spectrum consists of dominant eigenvalues or not. 
It zooms in on all the targets simultaneously, 
at practically the same rapid convergence rate.
The dimension of the 
subspaces, as well as the linear systems that need to be solved remained
unchanged throughout a fixed targeted domain $\domain$. Although the linear
systems are of the form $\phi_k B - A$, they are not shifts in the familiar sense. The
$\phi_k$s are not meant to be close to any eigenvalues but merely correspond to nodes
of a numerical quadrature rule. Under ideal situations, they are not near any
eigenvalues and none of the linear systems is ill-conditioned. Every one of these
features is distinct from those associated with popular non-Hermitian eigensolvers
such as unsymmetric Lanczos~\cite{parlett-taylor-liu-1985}, 
Arnoldi~\cite{lehoucq-sorensen-1996}, or 
Jacobi-Davidson~\cite{arbenz-hochstenbach-2004,sleijpen-vandervorst-2000}.
FEAST is fundamentally based on subspace iteration, 
whereas~\cite{sakurai-sugiura-2003,ikegami-sakurai-nagashima-2010},
despite their use of quadrature techniques, are more related to Krylov methods.
The quadratures there are used to approximate higher-order matrix moments.
In contrast, the quadratures used in FEAST are used to approximate the zeroth moments,
which correspond to spectral projectors.

\VS
The FEAST algorithms require the user to set a subspace dimension $p$, which should exceed $m$,
the number of eigenvalues in $\domain$. In practice, $p$ is often chosen based
on a priori knowledge or experience, or trial-and-error. A more elaborate theory exists,
similar to those detailed in~\cite{tang-polizzi-2013-FEAST} for the Hermitian case, on
estimation of the $m$.
For example, one can use
the eigenvalues of $\ctrans{\widehat{V}}B \widehat{U}$ 
($\widehat{U},\widehat{V}$ from Steps 4 and 5 of Bi-FEAST) 
to estimate the eigenvalue count $m$.

\VS
Opportunities for further work present themselves naturally, in the directions of
approximation theory, matrix analysis and parallel computing. At FEAST's core
is a rational function close to 1 inside a domain $\domain$, and
0 outside. Here we have used either a Gauss or trapezoidal quadrature rule to
construct this rational function. 
In general, possibility abounds for other quadrature rules,
either general or domain, $\domain$, specific
(see~\cite{bailey-borwein-2012} for example). Alternatively, one can view this as
a function approximation problem. 
Chebyshev polynomials~\cite{zhou-saad-2007,zhou-saad-tiago-chelikowsky-2006} 
which work well on the real line (for Hermitian problems) would not work on the
complex plane in terms of approximating the $\pi(\mu)$ function
in Equation (\ref{eqn:pi}):
Polynomials are analytic and must obey the
maximum modulus theorem (see~\cite{polya-latta-1974} for example). Rational approximation
can contribute fruitfully here. We have already seen one such case for Hermitian
problem where Zolotarev approximation is shown to outperform Gauss quadrature~\cite{viaud-thesis}.

\VS
In non-Hermitian matrix computations, it is customary to focus on the class of diagonalizable matrices.
How well a quadrature-based approximate spectral projector handles a general Jordan block,
and what the resulting implication on FEAST's convergence behavior in
the face of deficient eigenvectors will be,
are worthy pursuit that requires classical matrix and
perturbation analysis.

\VS
Last but not least, FEAST offers multiple levels of parallelism: 
multiple target domains, multiple linear systems, with multiple right hand sides.
Exploiting these parallelism fully, automatically, require much work still. On the highest
level, fast partitioning of a region in the complex plane to subregions, each containing
roughly the same number of eigenvalues, for the obvious sake of load balancing, 
is nontrivial. Challenging software engineering work is required to 
automatically distribute and
coordinate the linear solvers -- direct or iterative, sparse or dense -- on multiple right hand
sides, among multiple nodes, cores and threads.



\end{document}